\documentclass[letterpaper, 10 pt, journal, twoside]{ieeecolor}

\usepackage{lcsys}
\usepackage[left=0.7in,top=0.7in,right=0.7in,bottom=0.7in]{geometry}
\usepackage{mwe}
\usepackage{afterpage}
\usepackage{graphics} 
\usepackage{epsfig,url} 
\usepackage{times} 
\usepackage{bigints}

\usepackage{amsmath} 
\usepackage{amssymb,xcolor}  
\usepackage{mathtools,gensymb,amsfonts,amsthm,mathrsfs,comment,bm}
\usepackage{booktabs}
\usepackage[capitalise]{cleveref}

\theoremstyle{plain}
\newtheorem{theorem}{Theorem}
\newtheorem{lemma}[theorem]{Lemma}

\theoremstyle{definition}

\theoremstyle{remark}

\DeclarePairedDelimiter\abs{\lvert}{\rvert}
\DeclarePairedDelimiter\norm{\lVert}{\rVert}
\makeatletter
\let\oldabs\abs
\def\abs{\@ifstar{\oldabs}{\oldabs*}}
\let\oldnorm\norm
\def\norm{\@ifstar{\oldnorm}{\oldnorm*}}
\makeatother

\begin{document}
\newgeometry{left=0.70in,top=0.83in,right=0.70in,bottom=0.70in}
\afterpage{\aftergroup\restoregeometry}

\title{Boundary Output Feedback Stabilization for a Novel Magnetizable Piezoelectric Beam Model}
\author{Ahmet \"Ozkan \"Ozer \IEEEmembership{Member, IEEE}, Uthman Rasaq, Ibrahim Khalilullah
\thanks{A.\"O. \"Ozer is a faculty with Department of Mathematics, Western Kentucky University,
	Bowling Green, KY 42101, USA.  (e-mail: ozkan.ozer@wku.edu).}
\thanks{U. Rasaq and I. Khalilullah are graduate students at the Department of Mathematics, Western Kentucky University,
	Bowling Green, KY 42101, USA. (e-mails: uthman.rasaq114@topper.wku.edu, skmdibrahim.khalilullah504@topper.wku.edu).}
	\thanks{During A.O. Ozer's sabbatical at Friedrich-Alexander-Universität Erlangen-Nürnberg, we appreciate E. Zuazua's support. Additionally, we acknowledge the support from the U.S. National Science Foundation under Cooperative Agreement No. 1849213.}
}

\maketitle
\thispagestyle{empty}

\begin{abstract}
A magnetizable piezoelectric beam model, free at both ends, is considered.   Piezoelectric materials  have a strong interaction of electromagnetic and acoustic waves, whose wave propagation speeds differ substantially. The corresponding strongly-coupled PDE model describes the longitudinal vibrations and the total charge accumulation at the electrodes of the beam.  It is known that the PDE model  with appropriately chosen collocated state feedback controllers  is known  to have exponentially stable solutions. However, the collocated controller design is not always feasible since the performance of controllers  may not be good enough, and moreover, a small increment of feedback controller gains can easily make the closed-loop system unstable. Therefore, a non-collocated controller and observer design is considered for the first time for this model.  In particular, two state feedback controllers are designed at the right end to recover the states so that the boundary output feedback controllers can be designed as a replacement of the states with the estimate from the observers on the left end. By a carefully-constructed Lyapunov function, it is proved that the both the observer and  the observer error dynamics have  uniformly exponential stable solutions. This framework  offers a substantial foundation for the model reductions by Finite Differences. \end{abstract}


\section{Introduction}
{\color{black}

Piezoelectric materials, such as lead zirconate titanate (PZT), barium titanate, or lead titanate, are versatile smart materials that generate electric displacement directly proportional to applied mechanical stress \cite{Smith}. Primarily used in designing actuators, the drive frequency, a key component of the electrical input, dictates the speed of vibration or state changes in a piezoelectric beam \cite{Yang}. These materials also find application as sensors or energy harvesters, e.g. see \cite{Baur}.

While electrostatic approximations based on Maxwell's equations suffice for many piezoelectric applications, scenarios like piezoelectric acoustic wave devices highlight the dominance of magnetic effects. Existing models, relying solely on mechanical interactions, fail to capture the true vibrational dynamics in these devices. Hence, more accurate models incorporating electromagnetic interactions are essential, as demonstrated in \cite{D,T} and related references. Leveraging the complete set of Maxwell's equations for modeling ensures appropriate coupling of electromagnetic waves to mechanical vibrations, termed piezo-electro-magnetism by some researchers \cite{Yang}.

Introducing fully dynamic theory in modeling (acoustic) magnetizable piezoelectric beams significantly alters controllability dynamics compared to models derived through electrostatic approximation, known for exact controllability with a tip-velocity controller \cite{M-O, O-MCSS}. However, the fully dynamic model lacks controllability over high-frequency vibrational modes \cite{M-O}-\cite{AMOP}, \cite{Voss} using the same controller design. Consequently, a shift in the controller design is necessary to enable exact controllability.

\subsection{The PDE Model and Known  Stabilization Results}
Consider a magnetizable piezoelectric beam, clamped on one side and free to oscillate on the other, with sensors for tip velocity and tip current. The beam of length $L$ follows linear Euler-Bernoulli beam theory, with negligible transverse oscillations compared to longitudinal oscillations, manifesting as expansion and compression along the center line of the beam. Using $\alpha, \alpha_1, \gamma, \beta, \mu, \rho$ as positive piezoelectric material constants, define matrices
\begin{eqnarray}
\label{matrices1}
\begin{array}{ll}
 M:=\begin{bmatrix}
    \rho &0  \\
  0& \mu        \\
\end{bmatrix} , A:=\begin{bmatrix}
  \alpha       & -\gamma \beta   \\
   -\gamma\beta       & \beta  \\
\end{bmatrix}, ~\alpha=\alpha_1+\gamma^2\beta.
\end{array}
\end{eqnarray}
Here, $\alpha$ is influenced by piezoelectricity, distinct from $\alpha_1$ for fully elastic materials.  Noting time derivatives with dots and denoting longitudinal oscillations as $w(x,t)$ and total charge as $p(x,t)$, the motion is governed by a system of strongly-coupled partial differential equations, detailed in \cite{M-O},
\begin{eqnarray}
		\left\{
		\label{main}  \begin{array}{ll}
			\begingroup
M
\endgroup  \begin{bmatrix}
				{\ddot w}  \\
				{\ddot p}       \\
			\end{bmatrix}-
A \begin{bmatrix}
				w_{xx}  \\
				p_{xx}       \\
			\end{bmatrix}=0, ~x\in (0,L), t>0,\\
		A \begin{bmatrix}
				w_{x} \\
				p_{x}  \\
			\end{bmatrix}(0,t)=0,~~
A
 \begin{bmatrix}
				w_{x}  \\
				p_{x}   \\
			\end{bmatrix}(L,t) =\vec U(t), t>0,\\
\left(w,p,\dot w, \dot p \right)(x,0)=(w_0,p_0,w_1,p_1)(x),~x\in [0,L],\\
 \vec y(t)=\left(w,p,\dot w,\dot p \right)(0,t),~~t\ge 0.
		\end{array}\right.
	\end{eqnarray}
The constants $\rho, \alpha, \alpha_1, \beta, \gamma, \mu$ represent mass density, elastic and piezoelectric stiffness, beam impermeability, piezoelectric constant, and magnetic permeability, respectively. Controllers and measured output are denoted as $\vec U(t)$ and $\vec y(t)$.

The system with clamped-free boundary conditions lacks exponentially stability with a single controller, i.e., $\vec U(t)=-\begin{bmatrix}
k_w &0\\
0& k_p \\
\end{bmatrix}\begin{bmatrix}
\dot w \\
\dot p \\
\end{bmatrix}(L,t)$, where either $k_w=0$ and $k_p\ne 0$ or $k_w\ne 0$ and $k_p=0$. Thus, the necessity of two-state feedback controllers is emphasized \cite{M-O}. The associated closed-loop system imposes stringent stabilization conditions \cite{M-O}, with exponential stability achievable only for a limited subset of material parameters within a broad range \cite{O-MCSS}.

Recent studies highlight that employing two controllers ($k_w,k_p\ne 0$) results in exponential decay, supported by proving an observability inequality \cite{Wilson,Ramos} or utilizing a direct Lyapunov approach \cite{O-R}. In \cite{Ramos}, the proof employs a decomposition argument and an observability inequality with a suboptimal observation time, limiting the determination of optimal feedback gains for maximal decay rate. However, this limitation is overcome in \cite{O-R} by introducing a strategically-constructed Lyapunov function, providing explicit optimal feedback gains along with the maximal decay rate.

\subsection{Closed-loop Model Under Consideration}
 Define the solutions of \eqref{main} in the Hilbert space  $\mathcal H:=(H^1(0,L))^2\times (L^2(0,L))^2$. The new
 state feedback controller
\begin{eqnarray}
\begin{array}{ll}
\vec U(t)=\begin{bmatrix}
				-k_1 w -k_2 \dot w \\
				-k_3 p   -k_4 \dot p     \\
			\end{bmatrix}(L,t),  ~~k_1, k_2,k_3,   k_4>0,
\end{array}
\end{eqnarray}
and   the energy of the solutions
	\begin{eqnarray}\label{naturalenergy}
		\begin{array}{ll}
                     E(t)=&\frac{1}{2}\int^L_0 \left\{ M \begin{bmatrix}
				\dot w \\
				\dot p       \\
			\end{bmatrix} \cdot \begin{bmatrix}
				\dot w  \\
				\dot p       \\
			\end{bmatrix} + A  \begin{bmatrix}
				 w_x  \\
				p_x       \\
			\end{bmatrix} \cdot \begin{bmatrix}
				 w_x  \\
				p_x       \\
			\end{bmatrix} \right\} dx\\
		&	+ k_1 |w(L,t)|^2+ k_3 | p(L,t)|^2
		\end{array}
	\end{eqnarray}
make the closed-loop system  \eqref{main} dissipative.
 \begin{lemma}
The energy $E(t)$ in \eqref{naturalenergy} is dissipative, i.e.
\begin{eqnarray}\label{dissp}
\dot E(t)=-k_2 |\dot w(L,t)|^2-k_4|\dot p(L,t)|^2\le 0.
\end{eqnarray}
\end{lemma}
\begin{proof}
By taking the time derivative of  $E(t)$ along the solutions of \eqref{main} and by \eqref{matrices1},
\begin{eqnarray*}
\begin{array}{ll}
&\dot E(t)
=\int^L_0\left[\dot{w} \left(\alpha w_{xx}-\gamma\beta p_{xx}\right)+\dot{p}\left(\beta p_{xx}-\gamma\beta w_{xx}\right)\right]dx\\
&\quad  + \int^L_0\left[\alpha_1 w_x \dot {w_{x}}+\beta\left(\gamma w_x-p_x\right)\left(\gamma \dot{w_{x}}-\dot{p_{x}}\right)\right]dx \\
&\quad +k_1w(L,t)\dot{w}(L,t)+k_3p(L,t)\dot{p}(L,t)\\
&=\dot{w}(L)\left(\alpha w_x(L)- \gamma \beta p_x(L)\right)+\dot{p}(L)(\beta p_x(L) \\
&\quad-\gamma \beta w_x(L))+k_1w(L)\dot{w}(L)+k_2p(L)\dot{p}(L)
\end{array}
\end{eqnarray*}
and finally by the boundary conditions \eqref{main}, \eqref{dissp} follows.
\end{proof}
 Moreover, $E(t)$ in \eqref{naturalenergy} decays exponentially (which will be shown later in Theorem \ref{thm1}),  i.e. there exists constants $C,\sigma>0$ such that
$E(t)\leq C E(0)e^{-\sigma t}$ for all $t>0.$

Note that by the change of variable $x\to L-x$, the system \eqref{main} is equivalent to the following one
	\begin{eqnarray}
		\left\{
		\label{main-trans}  \begin{array}{ll}
M \begin{bmatrix}
				{\ddot w}  \\
				{\ddot p}       \\
			\end{bmatrix}-
A \begin{bmatrix}
				w_{xx}  \\
				p_{xx}       \\
			\end{bmatrix}=0, ~x\in (0,L), t>0,\\
		A \begin{bmatrix}
				w_{x} \\
				p_{x}  \\
			\end{bmatrix}(0,t)=\begin{bmatrix}
				k_1 w +k_2 \dot w \\
				k_3 p   +k_4 \dot p     \\
			\end{bmatrix}(0,t),\\
A
 \begin{bmatrix}
				w_{x}  \\
				p_{x}   \\
			\end{bmatrix}(L,t) =0, t>0,\\
\left(w,p,\dot w, \dot p \right)(x,0)=(w_0,p_0,w_1,p_1)(x),~x\in [0,L].
		\end{array}\right.
	\end{eqnarray}
By modifying the energy to the following one,
 	\begin{eqnarray}\label{modifiedenergy}
		\begin{array}{ll}
                   &  E(t)=\frac{1}{2}\int^L_0 \left\{ M \begin{bmatrix}
				\dot w \\
				\dot p       \\
			\end{bmatrix} \cdot \begin{bmatrix}
				\dot w  \\
				\dot p       \\
			\end{bmatrix} + A  \begin{bmatrix}
				 w_x  \\
				p_x       \\
			\end{bmatrix} \cdot \begin{bmatrix}
				 w_x  \\
				p_x       \\
			\end{bmatrix} \right\} dx\\
			&+ k_1 | w(0,t)|^2+ k_3 | p(0,t)|^2,
		\end{array}
	\end{eqnarray}
the closed-loop system \eqref{main-trans} is still dissipative.

\subsection{Non-collocated Observer Design and Error Tracking}
Define the following four positive-definite matrices
\begin{eqnarray}
\label{matrices2}
\begin{array}{ll}
 & K_{ij}=\begin{bmatrix}
    k_i &0  \\
  0& k_j        \\
\end{bmatrix}, ~~i=1,2,5,6, ~~j=i+2
\end{array}
\end{eqnarray}
where $k_1,\ldots, k_8>0$ are state feedback gains.
Since the controllers and observations are chosen non-collocated in \eqref{main}, like in \cite{Krstic,Ren},  designing stabilizing output feedback control requires designing
observers first as the following
\begin{eqnarray}\label{systemhat}
		\left\{
		 \begin{array}{ll}
M
  \begin{bmatrix}
				{\ddot {\hat{w}}}  \\
				{\ddot {\hat{p}}}       \\
			\end{bmatrix} -
A\begin{bmatrix}
				\hat{w}_{xx}  \\
				\hat{p}_{xx}       \\
			\end{bmatrix}=0,\\
\left(A \begin{bmatrix}
    \hat{w}_x\\
    \hat{p}_x
\end{bmatrix}= K_{13} \begin{bmatrix}
              {\dot {\hat{w}}}-{\dot {w}}\\
        {\dot {\hat{p}}} -   {\dot {p}}
    \end{bmatrix} +K_{24} \begin{bmatrix}
               {\hat{w}}- {w}\\
         \hat{p} -   p
    \end{bmatrix} \right)(0,t)\\
\left(A \begin{bmatrix}
    \hat{w}_x\\
    \hat{p}_x
\end{bmatrix}= -K_{57} \begin{bmatrix}
               {\dot {\hat{w}}}\\
        {\dot {\hat{p}}}
    \end{bmatrix}
     -K_{68} \begin{bmatrix}
               {\hat{w}}\\
         \hat{p}
    \end{bmatrix} \right)(L,t),\\
\left(\hat w, \hat p, \dot{\hat w}, \dot{\hat p} \right)(x,0)\\
\qquad\qquad =(w_0,p_0,w_1,p_1)(x),~x\in [0,L],
\end{array}\right.
\end{eqnarray}
with the energy of the observer defined by
\begin{eqnarray}\label{pseudoenergy1}
		\begin{array}{ll}
                     \hat E(t):=\frac{k_6}{2} \dot {\hat w}^2(L,t)+\frac{k_8}{2} \dot {\hat p}^2(L,t) &\\
			   \qquad+\frac{1}{2}\int^L_0 \left\{ M \begin{bmatrix}
				\dot {\hat w} \\
				\dot {\hat p}       \\
			\end{bmatrix} \cdot \begin{bmatrix}
				\dot {\hat w} \\
				\dot {\hat p}       \\
			\end{bmatrix} + A  \begin{bmatrix}
				 {\hat w}_x  \\
				{\hat p}_x       \\
			\end{bmatrix} \cdot \begin{bmatrix}
				 {\hat w}_x  \\
				{\hat p}_x       \\
			\end{bmatrix} \right\} dx.
		\end{array}
	\end{eqnarray}
Define the observer error for each variable  by
\begin{eqnarray}
\left\{
\begin{array}{ll}
e_1(x,t):=\hat{w}(x,t)-w(x,t)\\
e_2(x,t):=\hat{p}(x,t)-p(x,t)
\end{array}
\right.
\end{eqnarray}
where $e_1(x,t)$ and $e_2(x,t)$ satisfy the following system of error equations
\begin{eqnarray}
		\left\{
		\label{system-e}
		\begin{array}{ll} 
M \begin{bmatrix}
			\ddot e_1 \\
			\ddot e_2     \\
			\end{bmatrix} - A\begin{bmatrix}
				e_{1,xx}  \\
				e_{2,xx}       \\
			\end{bmatrix}=0,\\
 \left[A \begin{bmatrix}
    e_{1,x}\\
e_{2,x}
\end{bmatrix}= \begin{bmatrix}
              v_1\\
         v_2
    \end{bmatrix}:=K_{13} \begin{bmatrix}
             \dot  e_1\\
       \dot  e_2
    \end{bmatrix}   +K_{24} \begin{bmatrix}
               e_1\\
      e_2
    \end{bmatrix} \right](0,t),\\
 A \begin{bmatrix}
    e_{1,x}\\
e_{2,x}
\end{bmatrix}(L,t)= 0,\\
\left(e_1, e_2, \dot e_1,\dot e_2 \right)(x,0)=0,~x\in [0,L],
\end{array}\right.
\end{eqnarray}
with the energy of the observer error defined by
\begin{eqnarray}\label{pseudoenergy2}
		\begin{array}{ll}
                     E_e(t):=\frac{k_2}{2}e^2_1(0,t) +\frac{k_4}{2}e^2_2(0,t)\\
                     +\frac{1}{2}\int^L_0 \left\{ M \begin{bmatrix}
				\dot e_1\\
				\dot e_2    \\
			\end{bmatrix} \cdot \begin{bmatrix}
				\dot e_1\\
				\dot e_2    \\
			\end{bmatrix} + A  \begin{bmatrix}
				 e_{1,x}  \\
				 e_{2,x}      \\
			\end{bmatrix} \cdot \begin{bmatrix}
				 e_{1,x}  \\
				 e_{2,x}      \\
			\end{bmatrix} \right\} dx. &
		\end{array}
	\end{eqnarray}
Note that the system \eqref{system-e} is the same as the system \eqref{main-trans}, which is also equivalent to	\eqref{main}.	
		\subsection{Our Contributions}
To the best of our knowledge, this work marks the inaugural effort towards a non-collocated observer design for a strongly-coupled and fully-dynamic (magnetizable) piezoelectric beam model. Through the design of tip velocity and total current observers at the left end of the beam,  we establish exponential stability in overall vibrations towards the equilibrium state. The proof of our main result relies on a meticulously chosen Lyapunov function. Importantly, our proof technique obviates the need for spectral analysis, enabling the design of an order-reduced model reduction through Finite Differences, reminiscent of the approach in \cite{Ren}. This model reduction stands as a primary objective of the paper.

}
\section{Exponential Stability by Lyapunov Approach}
The following technical lemmas are needed for the exponential stability result.
Let $E(t):=\hat E(t)+ E_e(t).$
\begin{lemma}
With $v_1(0,t)$ and $v_2(0,t)$, defined by \eqref{system-e}, the energies $\hat E(t)$ and $E_e(t)$ in \eqref{pseudoenergy1} and \eqref{pseudoenergy2}  satisfy
\begin{eqnarray}\label{dissip}
    \begin{array}{ll}
     \dot {\hat E}(t) = -k_5|\dot{\hat w}(0,t)|^2-k_7|\dot{\hat p}(0,t)|^2&\\
     \qquad\quad~ -\dot{\hat w}(0,t)  v_1(0,t)-\dot{\hat p}(0,t)  v_2(0,t),&\\
        \dot {E_e}(t) =-k_1 |\dot e_1(0,t)|^2-k_3 |\dot e_2(0,t)|^2\le 0.
    \end{array}
\end{eqnarray}

\end{lemma}
\begin{proof} Take the derivative of $\hat{E}(t)$ along with the solutions of \eqref{systemhat} together with \eqref{matrices1}
\begin{eqnarray*}
   \begin{array}{ll}
  \dot {\hat E}(t)= \int^L_0 \left[\rho \dot{ \hat w}_{1} \ddot{ \hat w}_{1}+\mu \dot{\hat w}_{2} \ddot{\hat w}_{2} +\alpha_1  \hat{w}_{1,x} \dot{\hat w}_{1,x}\right. \\
     \qquad\qquad \left. +\beta \left(\gamma \hat{w}_{1,x}-\hat{w}_{2,x}\right)\left(\gamma\hat{w}_{1,xt}-\hat{w}_{2,xt}\right)\right]dx\\
     \qquad\qquad +k_6 \hat{w}_{1,t}(L) \hat{w}_1(L)+k_8\hat{w}_{2,t}(L)\hat{w}_2(L)\\
     =\int^L_0 \left[\hat{w}_{1,t} \left(\alpha \hat{w}_{xx}  -\gamma \beta \hat{p}_{xx}\right) +\beta\hat{w}_{2,t}\left(\hat{p}_xx-\gamma  \hat{w}_{xx}\right)\right. \\
     \quad \left.+ \alpha_1 \hat{w}_{1,x}\hat{w}_{1,xt} + \beta\left(\gamma  \hat{w}_{1,x}-\hat{w}_{2,x}\right)(\gamma \hat{w}_{1,xt} -\hat{w}_{2,xt}) \right]dx\\
    \quad + k_6\hat{w}{1,t}(L)\hat{w}_1(L)+k_8\hat{w}_{2,t}(L)\hat{w}_2(L),
    \end{array}
\end{eqnarray*}
and the first equality  in \eqref{dissip} follows from the boundary conditions \eqref{systemhat}.  Next, take  the derivative of ${E}_e(t)$ along the solutions of \eqref{system-e}  with the boundary conditions and \eqref{matrices1} to get the second equality in \eqref{dissip},
\begin{eqnarray*}
    \begin{array}{ll}
       &  \dot {E_e}(t) =\int^L_0 \big(\rho \dot e_{1} \ddot e_{1}+\mu \dot e_{2} \ddot e_{2}+ 2\alpha_1 e_{1,x} \dot e_{1,x}\\
     &\qquad\qquad +\beta \big(\gamma \dot e_{1,x}-\dot e_{2,x}\big)\big(\gamma e_{1,x}-e_{2,x}\big)\big)dx\\
     &\qquad +k_2e_1(0,t) \dot e_{1}(0,t)+k_4e_2(0,t) \dot e_{2}(0,t) \\
     &=\left\{k_2e_1 \dot e_{1}-\dot e_{1}\left(\alpha e_{1,x}-\gamma \beta e_{2,x}\right))\right. \\
     &\quad \left.+k_4e_2 \dot e_{2}- \dot e_{2} (\beta e_{2,x} -\gamma \beta e_{1,x})\right\}(0,t). \qedhere
     \end{array}
\end{eqnarray*}
\end{proof}

For $C_e, \epsilon_1, \epsilon_2, N_1, N_2>0,$ define the  Lyapunov function $L(t):=C_e  L_e(t)+ \hat L(t)$ where
\begin{eqnarray}
\label{Lyap}&&\left\{\begin{array}{ll}
L_e(t):=E_e(t) +\epsilon_1\left\{\psi_{1,1}(t)+\psi_{1,2}(t) \right.&\\
\qquad \qquad \qquad\qquad \left. +N_{1}(\psi_{2,1}(t)+\psi_{2,2}(t))\right\},\\
\hat L(t):= \hat{E}(t)+\epsilon_2\left\{ \psi_{3,1}(t)+\psi_{3.2}(t) \right.&\\
\qquad \qquad\qquad \qquad\left. +N_{2}(\psi_{4,1}(t)+\psi_{4,2}(t)) \right\},\\
\end{array}\right.
\end{eqnarray}
\begin{eqnarray}
\label{19}\left\{\begin{array}{ll}
\psi_{1,1}(t):=\rho \int^L_0 (x-L)\dot e_{1} e_{1,x}~dx, \\
\psi_{1,2}(t):=\mu \int^L_0(x-L)\dot e_{2} e_{2,x}~dx,\\
\psi_{2,1}(t):=\rho e_1(0,t)\int^L_0 \dot e_1~dx+\frac{k_1}{2}e_1(0,t)^2,\\
\psi_{2,2}(t):=\mu e_2(0,t)\int^L_0 \dot e_2~dx+\frac{k_3}{2}e_2(0,t)^2,\\
\psi_{3,1}(t):=\rho \int^L_0 \left(x+2L\right)  \dot{\hat{w}} \hat{w}_x~dx\\
\qquad-\frac{3L}{\alpha_1}k_5k_6\hat{w}(L,t)^2 -\frac{3L}{\alpha_1}\gamma^2 k_7k_8\hat{p}(L,t)^2,\\
\psi_{3,2}(t):=\mu \int^L_0  (x+2L)\dot{\hat{p}} ~\hat{p}_{x}dx-\frac{3L}{2\beta}k_7k_8 \hat{p}(L,t)^2,\\
\psi_{4,1}(t):=\hat{w}(L,t)\int^L_0\dot {\hat{w}}~dx+\frac{k_5}{2} \dot {\hat{w}}(L,t)^2,\\
\psi_{4,2}(t):=\hat{p}(L,t)\int^L_0 \dot{\hat{p}}~dx+\frac{k_7}{2} \dot{\hat{p}}(L,t)^2.
\end{array}\right.
\end{eqnarray}

\begin{lemma}  \label{lemm2}
The Lyapunov function $L(t)$ in \eqref{Lyap} and the energy $E(t)$ in \eqref{naturalenergy} are equivalent, i.e. there exists $p_1,p_2>0$ such that
\begin{equation}\label{lemres-c}
    \begin{array}{ll}
         p_1    E (t)\le L (t) \le p_2 E(t),
    \end{array}
\end{equation}
where for all $0 < \epsilon_1<C^{-1}_1 ,0< \epsilon_2<C^{-1}_2,$ and
\begin{equation}\label{c1}
    \begin{array}{ll}
       p_1 := \min \{C_e(1-\epsilon_1 C_1),1-\epsilon_2 C_2\},\\
     p_2 := \max \{C_e(1+\epsilon_1 C_1),1+\epsilon_2 C_2\},\\
       C_1 := L \max \left\{{ \sqrt{\frac{\rho}{\alpha_1}}+\sqrt{\frac{\mu \gamma^2}{\alpha_1}}, \sqrt{\frac{\mu}{\beta}}+\sqrt{\frac{\mu \gamma^2}{\alpha_1}}}\right\}  \\
       ~~ +N_1\max \left\{ \frac{1+k_1}{k_2},\frac{1+k_3}{k_4},L\right\},\\
       C_2 :=(3L) \max\left\{\sqrt{\frac{\rho}{\alpha_1}}+\sqrt{\frac{\mu \gamma^2}{\alpha_1}},\sqrt{\frac{\mu}{\beta}}+\sqrt{\frac{\mu \gamma^2}{\alpha_1}}\right.\\
       ~~ \left.\big(\frac{2 \gamma^2}{\alpha_1}+\frac{1}{\beta}\big),\frac{2}{\alpha_1}\right\} + N_2 \max \left\{\frac{\frac{1}{\rho}+k_5}{k_6},\frac{\frac{1}{\mu}+k_7}{k_8},L\right\}.
    \end{array}
\end{equation}
\end{lemma}
\begin{proof}
By  Cauchy-Schwartz's, H\"older's, Young's and Minkowski's inequalities,
\begin{eqnarray}\label{psi_1}
    \begin{array}{ll}
      &|\psi_{1,1}(t)+\psi_{1,2}(t)|\le \frac{L}{2} \sqrt{\frac{\mu}{\beta}} \left. \mu \int^L_0 |e_{2,t}|^2  \right.\\
      &\qquad\qquad \left. + \beta \int^L_0 (e_{2,x}- \gamma e_{1,x})^2\right) dx\\
      & \quad+\frac{L}{2} \sqrt{\frac{\mu \gamma^2}{\alpha_1}} \big(\mu \int^L_0 |e_{2,t}|^2dx+ \int^L_0 \alpha_1 |e_{1,x}|^2 dx \big)\big) \\
      &\le L \max \left\{\sqrt{\frac{\rho}{\alpha_1}}+\sqrt{\frac{\mu \gamma^2} {\alpha_1}}, \sqrt{\frac{\mu}{\beta}}+\sqrt{\frac{\mu \gamma^2}{\alpha_1}}\right\} E_e(t).
    \end{array}
\end{eqnarray}
Next, by Young's and Cauchy-Schwartz's inequalities,
\begin{eqnarray}\label{psi_2}
    \begin{array}{ll}
    &|\psi_{2,1}(t)+\psi_{2,2}(t)|\le\frac{\rho+k_1}{2k_2}k_2 e^2_1(0)\\
    &\quad +\frac{\mu+k_3}{2k_4}k_4 e^2_2(0) + \frac{L}{2} \big(\int^L_0 \rho {\dot e_1}^2 dx + \int^L_0 \mu {\dot e_2}^2 dx \big) \\
    & \le \max \left\{\frac{\rho+k_1}{k_2},\frac{\mu+k_3}{k_4},L \right\} E_e(t).
    \end{array}
\end{eqnarray}
Finally,  \eqref{psi_1} together with \eqref{psi_2} lead to
\begin{eqnarray}\label{le}
\begin{array}{ll}
 |L_e(t)| \le E_e(t) + \epsilon_1 |\psi_{1,1}+\psi_{1,2}|+ \epsilon_1 N_1|\psi_{2,1}+\psi_{2,2}| \\
 \le \left(1+\epsilon_1 \left( L\max \left\{ \sqrt{\frac{\rho}{\alpha_1}}+\sqrt{\frac{\mu \gamma^2}{\alpha_1}}, \sqrt{\frac{\mu}{\beta}}+\sqrt{\frac{\mu \gamma^2}{\alpha_1}}\right) \right. \right.\\
 \left.\qquad +N_1\max \left\{ \frac{\rho+k_1}{k_2},\frac{\mu+k_3}{k_4},L\right\} \right)E_e(t).
\end{array}
\end{eqnarray}
Therefore,
\begin{equation}\label{lemres-c1}
    \begin{array}{ll}
     (1-\epsilon_1 C_1)    E_e (t)\le L_e (t) \le (1+ \epsilon_1 C_1) E_e (t).
    \end{array}
\end{equation}
By the same set of inequalities used above, the following estimates can be obtained
\begin{eqnarray}\label{23}
    \begin{array}{ll}
    |\psi_{3,1}(t)+\psi_{3,2}(t)|\\
   \le \frac{3L}{2} \left\{ \frac{k_7 k_8}{2} \hat{p}^2(L) + \left[ \int^L_0 \sqrt{\frac{\mu \gamma^2}{\alpha_1}} \left(  \mu |{\dot {\hat p}}|^2+ \alpha_1 |\hat{p}_x|^2  \right)\right.\right.\\
    \quad\left. \left.+ \sqrt{\frac{\mu}{\beta}} \left( \mu|{\dot {\hat p}}|^2 + \beta  (\hat{p}_x -\gamma \hat{w}_x )^2\right)\right] dx \right\}\\
    \le 3L \max \left\{\sqrt{\frac{\rho}{\alpha_1}}+\sqrt{\frac{\mu \gamma^2}{\alpha_1}}, \sqrt{\frac{\mu}{\beta}}+\sqrt{\frac{\mu \gamma^2}{\alpha_1}}\right.,\\
    \left. \qquad\qquad\qquad\qquad\big(\frac{2 \gamma^2}{\alpha_1}+\frac{1}{\beta}\big),\frac{2}{\alpha_1}\right\} \hat{E}(t)\\
    \end{array}
\end{eqnarray}
 and, analogously,
\begin{eqnarray}\label{24}
    \begin{array}{ll}
|\psi_{4,1}(t) +\psi_{4,2}(t)|
\le \max \left\{\left(\frac{\rho+k_5}{2k_6}\right)k_6 {\hat{w}}^2(L,t)\right. \\
\left.\qquad+\left(\frac{\mu+k_7}{2k_8}\right)k_8 {\hat{p}}^2(L,t) +\frac{L}{2} \int^L_0 (\rho \dot {\hat w}^2 + \mu \dot{\hat p}^2 \big)dx \right\}\\
\le \max \left\{\frac{\rho+k_5}{k_6},\frac{\mu+k_7}{k_8},L \right\}\hat{E}(t).
     \end{array}
\end{eqnarray}
Considering \eqref{23} together with \eqref{24} yields
\begin{eqnarray}\label{lhat}
\begin{array}{ll}
& \hat{L}(t) \le \hat{E}(t)+ \epsilon_2(\psi_{3,1}+\psi_{3,2})+ \epsilon_2 N_2(\psi_{4,1} +\psi_{4,2})\\
 &  ~~ \le \hat{E}(t)+ \epsilon_2\left[ N_2 \max \left\{\frac{\rho+k_5}{k_6},\frac{\mu+k_7}{k_8},L\right\}\right.\\
   &  \qquad +3L\max\left\{\sqrt{\frac{\rho}{\alpha_1}}+\sqrt{\frac{\mu \gamma^2}{\alpha_1}}, \sqrt{\frac{\mu}{\beta}}+\sqrt{\frac{\mu \gamma^2}{\alpha_1}},\right.\\
    & \qquad  \qquad \qquad \qquad  \qquad \left.\left.(\frac{2\gamma^2}{\alpha_1}+\frac{1}{\beta}),\frac{2}{\alpha_1}\right\}\right]\hat{E}(t)\\
    & ~~\le \left( 1+ \epsilon_2\big[ N_2 \max \left\{\frac{\rho+k_5}{k_6},\frac{\mu+k_7}{k_8},L\right\} \right.\\
     &\qquad +3L\max\left\{\sqrt{\frac{\rho}{\alpha_1}}+\sqrt{\frac{\mu \gamma^2}{\alpha_1}}, \sqrt{\frac{\mu}{\beta}}+\sqrt{\frac{\mu \gamma^2}{\alpha_1}},\right.\\
     &\left.\left. \left.\qquad\qquad\qquad\qquad\qquad (\frac{2\gamma^2}{\alpha_1}+\frac{1}{\beta}),\frac{2}{\alpha_1}\right\}\right]\right) \hat{E}(t),
         \end{array}
\end{eqnarray}
and therefore,
\begin{equation}\label{lemres-c2}
    \begin{array}{ll}
        (1-\epsilon_2 C_2)    \hat{E}(t)\le \hat{L}(t) \le (1+ \epsilon_2 C_2) \hat{E}(t).
    \end{array}
\end{equation}
Finally \eqref{lemres-c} follows from \eqref{lemres-c1} and \eqref{lemres-c2}.
\end{proof}
\begin{lemma}  \label{lemm4}
Let $0<L<2,$ $ \delta_1<2\epsilon_2 \rho,$  $\delta_2<2\epsilon_2 \mu,$ and

\begin{equation*}\label{caps1}
    \begin{array}{ll}
  C_e> \max \left[\frac{k_1^2\left|\frac{1}{\delta_1}+\epsilon_2 N^2_2 +\frac{\gamma^2 \epsilon_2}{\alpha_1}-\frac{2\epsilon_2}{\alpha_1}\right|}{k_1-\epsilon_1\left[\left[\frac{\color{black}\rho}{2}+\frac{2k^2_1}{\alpha_1}\right]L+N^2_1 \rho\right]},  \frac{k^2_2\left|\frac{1}{\delta_1}+ \epsilon_2N^2_2+\frac{\gamma^2 \epsilon_2}{\alpha_1}-\frac{2 \epsilon_2}{\alpha_1}\right|}{\epsilon_1 \left[N_1k_2-\left[\frac{k_2}{2}+\frac{Lk^2_2}{\alpha_1}\right]\right]}\right],\\
 ~~\frac{k_3^2\left|\frac{1}{\delta_2}+\epsilon_2 N^2_2 +\frac{2\epsilon_2 \gamma^2}{\alpha_1}-\frac{2\epsilon_2}{\beta}\right|}{k_3-\epsilon_1 \left[\left(\frac{\color{black}\mu}{2}+\frac{\alpha}{\alpha_1 \beta}k^2_3\right)L+ \mu N^2_1\right]},\left. \frac{k^2_4\left|\frac{1}{\delta_2}+\frac{2 \epsilon_2 \gamma^2}{\alpha_1}+\epsilon_2 N^2_2 -\frac{2 \epsilon_2}{\beta}\right|}{\epsilon_1 \left(N_1k_4-\left(\frac{k_4}{2}+\frac{\alpha L}{\alpha_1\beta}k^2_4\right)\right)} \right],\\
  \epsilon_1<\min\left[\frac{k_1}{\left(\frac{\rho}{2}+\frac{2k^2_1}{\alpha_1}\right)L+\rho N^2_1}, \frac{k_3}{\left(\frac{\mu}{2}+\left(\frac{\alpha_1+2\gamma^2 \beta}{\alpha_1 \beta}\right)k^2_3\right)L+\mu N^2_1}\right],\\
  \epsilon_2<\min \left[\frac{\frac{2}{3k_7}}{L\left[\rho(1+N_1^2)+\frac{2k_5^2}{\alpha_1}\right]}, \frac{\frac{2}{3k_7}}{L\left[\mu(1+N_2^2)+\left(\frac{2\gamma^2}{\alpha_1}+\frac{1}{\beta}\right)k_7^2\right]}\right],\\
  N_1> \frac{1}{2}+L\max\left[\frac{k_2}{\alpha_1}, \frac{\alpha k_4}{\alpha_1\beta} \right],\\
    N_2>\frac{1}{2}+\max \left[\frac{1}{2k_6}+ \frac{3L k_6}{\alpha_1},\frac{1}{2k_8}+ \frac{3L\left(\alpha+\gamma^2\beta\right) k_8}{2\alpha_1\beta} \right].
  \end{array}
\end{equation*}

Then, the Lyapunov function $\hat L(t)$  satisfies
\begin{equation}\label{Ldot}
    \begin{array}{ll}
    \dot {L}(t)\le  C_e  \epsilon_1 \left(1-\frac{L}{2}\right)E_e(t)-\frac{2\epsilon_2}{3} \hat E(t).
    \end{array}
\end{equation}
\end{lemma}
 \begin{proof} First, we estimate the time derivative of each $\psi-$function. By the integration by parts along with inequality $(a+b)^2 \leq 2(a^2 + b^2)$,
\begin{eqnarray}\label{psidot1}
    \begin{array}{ll}
  |\dot \psi_{1,1}(t)+\dot\psi_{1,2}(t)|\le   L\left[\left(\frac{\rho}{2}+\frac{2k_1^2}{\alpha_1} \right)  |\dot e_{1}(0,t)|^2\right. \\
\left.+\left(\frac{\mu}{2}+ \frac{\alpha+\gamma^2 \beta}{\alpha_1\beta} k_3^2\right) |\dot e_{2}(0,t)|^2 \right]   \\
+\left(\frac{k_2}{2}+\frac{2Lk_2^2}{\alpha_1} \right)  |e_{1}(0,t)|^2  \\
 +\left(\frac{k_4}{2}+ \frac{\alpha+\gamma^2 \beta}{\alpha_1\beta }L k_4^2\right) |e_{2}(0,t)|^2  -E_e(t)
    \end{array}
\end{eqnarray}
Next, by the Cauchy-Schwartz's and Young's inequalities,
\begin{eqnarray}\label{psidot2}
\begin{array}{ll}
  & N_1|\dot\psi_{2,1}(t)+\dot\psi_{2,2}(t)| \le N_1 \dot e_1(0) \int^L_0 \rho \dot e_1 dx\\
  &\quad + N_1 \dot e_2(0) \int^L_0 \mu \dot e_2 dx\\
  & \quad +N_1 \dot e_1(0) \int^L_0 \rho \big(\alpha e_{1,xx}-\gamma \beta e_{2,xx}\big)dx\\
  & \quad +N_1 \dot e_2(0) \int^L_0 \big(\beta e_{2,xx}-\gamma \beta e_{1,xx}\big)dx \\
  & \quad +N_1 k_1 e_1(0) \dot e_1(0,t) + N_1 k_3 e_2(0) \dot e_2(0,t) \\
  & \le  N^2_1 \bigg(\rho \dot e_1^2(0,t)+ \mu \dot e_2^2(0,t)\bigg)\\
  &\quad +\frac{L}{4} \int^L_0 \big(\rho \dot e^2_1+ \mu \dot e^2_2 \big)dx  \\
  & \quad-N_1 \big(k_2 e^2_1(0,t) +k_4 e^2_2(0,t) \big).
    \end{array}
\end{eqnarray}
Putting \eqref{psidot1} and \eqref{psidot2} together  $|\dot L_e(t)|$ can be estimates as the following
\begin{eqnarray*}
    \begin{array}{ll}
         &|\dot L_e(t)| \le \dot E_e(t) + \epsilon_1 \left[|\dot \psi_{1,1}+\dot \psi_{1,2}|+N_1|\dot \psi_{2,1}+\dot \psi_{2,2}|\right] \\
         & \le -k_1 {\dot e_1}^2(0,t)-k_3 {\dot e_2}^2(0,t)+\epsilon_1 \bigg[\big(\frac{\rho}{2}+
         \frac{2k^2_1}{\alpha_1}\big)L {\dot e_1}^2(0,t)\\
         &\quad+\left(\frac{\mu}{2}+\big(\frac{\alpha_1+2 \gamma^2 \beta}{\alpha_1}\big)k^2_3\right)L {\dot e}^2_2(0,t)\\
         &\quad +\left(\frac{k_4}{2}+\big(\frac{\alpha_1+2 \gamma^2 \beta}{\alpha_1 \beta}\big)k^2_4 L \right)e^2_2(0,t)-E_e(t)\\
         &\quad +N^2_1(\rho {\dot e_1}(0,t)+\mu {\dot e_2}^2(0,t))+\frac{L}{4} \int^L_0 (\rho {\dot e_1}^2+\mu {\dot e_2}^2)dx \\
         &\quad -N_1\big(k_2 e^2_1(0,t)+k_4 e^2_2(0,t)\big)\bigg].
            \end{array}
\end{eqnarray*}
Finally, by factoring out the coefficients of the energy terms, we get the following
\begin{eqnarray}
\label{Ledot}
    \begin{array}{ll}
         |\dot L_e(t)|  \le-\epsilon_1\big(1-\frac{L}{2}\big)E_e(t)\\
         ~-\big[k_1-\epsilon_1 \big(\frac{\rho}{2}+\frac{2k^2_1}{\alpha_1}\big)L -\epsilon_1 N^2_1 \rho \big]\dot e^2_1 (0,t)\\
         ~ -\big[k_3-\epsilon_1 \left[\frac{\mu}{2}+\left[\frac{\alpha_1+2\gamma^2 \beta}{\alpha_1 \beta})k^2_3 \right]L-\epsilon_1 \mu N^2_1 \right]\dot {e}^2_2 (0,t)\\
         ~-\epsilon_1\big[N_1 k_2-\big(\frac{k_2}{2}+\frac{2LK^2_2}{\alpha_1}\big)\big]e^2_1(0,t) \\
         ~ -\epsilon_1 \big[N_1 k_4-\big(\frac{k_4}{2}+\big(\frac{2 \gamma^2}{\alpha_1}+\frac{1}{\beta}\big)k^2_4 L\big)\big]e^2_2(0,t).
    \end{array}
\end{eqnarray}
Next, take the time derivative and  integrate by
parts together so that $\dot{\hat L}$ is estimated as follows
\begin{eqnarray}\label{psidot3}
    \begin{array}{ll}
  |\dot\psi_{3,1}(t)+\dot\psi_{3,2}(t)|\le  -E_w(t)+\rho {\dot {\hat w}}^2(0,t) -\mu {\dot {\hat p}}^2(0,t)\\
    \quad - \frac{1}{\alpha_1}\big(v_1(0,t)+ \gamma v_2(0,t) \big)^2- \frac{1}{\beta}v^2_2(0,t)\\
  \quad+ \frac{3L}{2}\left[\big(\rho+\frac{2k^2_3}{\alpha_1}\big) {\dot{\hat w}}^2(L,t) \right.\\
  \quad\left.\quad+\left(\mu+\left(\frac{2 \gamma^2}{\alpha_1}+\frac{1}{\beta}\right)k^2_7\right){\dot {\hat p}}^2(L,t)\right]\\
  \quad +3L\left[\left(\big(\frac{k^2_6}{\alpha_1}+\frac{k_6}{6L} )\hat{w}^2(L,t) \right.\right.\\
  \quad \left.\left.\quad+\big(\frac{\gamma^2}{\alpha_1}+\frac{1}{2\beta}\right)k^2_8+ \frac{k_8}{6L}\big) \hat{p}^2(L,t)\right]\\
    \end{array}
\end{eqnarray}
where the  generalized Cauchy-Schwartz's inequality $ab \leq \frac{a^2}{2\epsilon}+\frac{\epsilon b^2}{2}$ for $\epsilon = 3L$ and for  $\epsilon = 1 $ are used successively. Analogously,
\begin{eqnarray}\label{psidot4}
    \begin{array}{ll}
  &N_2|\dot\psi_{4,1}(t)+\dot\psi_{4,2}(t)|\le \frac{1}{6} \int^L_0 (\rho |{\dot {\hat w}}|^2 +\mu |{\dot {\hat p}}|^2)dx\\
    &\quad +\left[ -N_2\left(k_6 \hat{w}^2 +k_8 \hat{p}^2\right)+\frac{3LN^2_2 }{2}\left(\rho {\dot {\hat w}}^2+\mu {\dot {\hat p}}^2\right)\right.\\
  & \left. \quad + \frac{1}{2} \left(\hat{w}^2 +\hat{p}^2\right)\right](L,t)+ \frac{N^2_2}{2}\left(v^2_1 (0,t)+ v^2_2 (0,t)\right).
    \end{array}
\end{eqnarray}
Using \eqref{psidot3} and \eqref{psidot4} together with the generalized Cauchy-Schwartz's inequality for the terms $\dot{\hat w}(0)v_1(0)$ and $\dot{\hat p}(0)v_2(0)$ with $\epsilon=\delta_1$ and $\delta_2$, $\dot{\hat L}(t),$ respectively,
\begin{eqnarray*}\label{hatLdot}
    \begin{array}{ll}
& \dot{\hat L}(t) \le \dot{\hat E}(t)+ \epsilon_2[|\dot\psi_{3,1}+\dot\psi_{3,2}|+|N_2(\dot\psi_{4,1}+\dot\psi_{4,2})|]\\
& \le -\frac{\epsilon_2}{2} \int_0^L \left[ \frac{2}{3} \left(\rho\dot{\hat w}^2+ \mu\dot{\hat p}^2\right)+\alpha_1 \hat{w}^2_x +\beta (\gamma \hat{w}_x-\hat{p}_x)^2 \right]dx\\
&~-\left[k_5-\frac{3L\epsilon_2}{2}(\rho+\frac{2k^2_5}{\alpha_1}+N^2_2) \rho\right]\dot{\hat w}^2(L,t)\\
&~-\left[k_7-\frac{3L\epsilon_2}{2}\big(\mu+(\frac{2\gamma^2}{\alpha_1}+\frac{1}{\beta})k^2_7+N^2_2) \mu\right]\dot {\hat p}^2(L,t)\\
&~-\epsilon_2\left[N_2 k_6-\frac{3Lk^2_6}{\alpha_1}-\frac{k_6}{2}-\frac{1}{2}\right]\hat{w}^2(L,t)\\
&~ -\epsilon_2\left[N_2 k_8-(\frac{3L\gamma^2}{\alpha_1}+\frac{1}{2\beta})k^2_8-\frac{1}{2} \right]\hat{p}^2(L,t)\\
&~ -\left[\epsilon_2 \rho-\frac{\delta_1}{2}\right]\dot{ \hat w}^2(0)-\left[\epsilon_2 \mu -\frac{\delta_2}{2}\right]\dot{ \hat p}^2(0,t)\\
&~+\left[\frac{1}{2\delta_1}-\frac{\epsilon_2}{\alpha_1}+\frac{\epsilon_2 N^2_2}{2}+\frac{\gamma^2 \epsilon_2}{2 \alpha_1}\right]v^2_1(0,t)\\
&~+\left[\frac{1}{2 \delta_2}-\frac{\epsilon_2 \gamma^2}{\alpha_1}+\frac{N^2_2 \epsilon_2}{2}-\frac{\epsilon_2}{\beta}+\frac{2\epsilon_2}{\alpha_1}\right]v^2_2(0,t).
    \end{array}
\end{eqnarray*}
This, together with  \eqref{Ledot} and $v_1(0,t) = k_1 \dot{e}_1(0,t) + k_2{e}_1 (0,t),$  $v_2(0,t) = k_3 \dot{e}_2 (0,t) + k_4{e}_2 (0,t)$  lead to
\begin{eqnarray*}\label{Ldot1}
    \begin{array}{ll}
&|\dot {L}(t)| \le -C_e \epsilon_1\left(1-\frac{L}{2}\right)E_e(t)-\frac{2\epsilon_2}{3}E_{\hat{w}}(t)\\
&\quad-\left[ C_e\left(k_1- \epsilon_1 \big(\frac{\rho}{2}+\frac{2k^2_1}{\alpha_1}\big)L-\epsilon_1 N^2_1 \rho \right) \right.\\
&\qquad \left. -k^2_1\left(\frac{1}{\delta_1}-\frac{2\epsilon_2}{\alpha_1}+\epsilon_2N^2_2+\frac{\gamma^2 \epsilon_2}{\alpha}\right) \right]{\dot e_1}^2(0,t)\\
&\quad-\left[C_e\big(k_3-\epsilon_1\left(\frac{\mu}{2}+\big(\frac{\alpha_1+2\gamma^2 \beta}{\alpha_1 \beta}\big)k^2_3 \big)L-\epsilon_1 \mu N^2_1\right)\right.\\
&\qquad \left. -k^2_3\big(\frac{1}{\delta_2}+\frac{2\epsilon_2 \gamma^2}{\alpha_1}+\epsilon_2 N^2_2-\frac{2\epsilon_2}{\beta}\big)\right]{\dot e_2}^2(0,t)\\
&\quad-\left[ C_e \epsilon_1 \big(N_1k_2-\big(\frac{k_2}{2}+\frac{2Lk^2_2}{\alpha_1}\big)\big)\right.\\
&\quad \left. \qquad -k^2_2\big(\frac{1}{\delta_1}-\frac{2 \epsilon_2}{\alpha_1}+ \epsilon_2N^2_2+\frac{\gamma^2 \epsilon_2}{\alpha_1}\big)\right]e^2_1(0,t)\\
&\quad-\left[C_e \epsilon_1 \big(N_1k_4-\left(\frac{k_4}{2}+\big(\frac{2 \gamma^2}{\alpha_1}+\frac{1}{\beta}\big)k^2_4L\right)\right. \\
&\quad-k^2_4\big(\frac{1}{\delta_2}+\frac{2 \epsilon_2 \gamma^2}{\alpha_1}+N^2_2 \epsilon_2-\frac{2 \epsilon_2}{\beta}\big)\big]e^2_2(0,t)\\
&\quad-\epsilon_2\big[N_2 k_6-\frac{k_6}{2}-\frac{ 3Lk^2_6}{\alpha_1}-\frac{k_6}{2}-\frac{1}{2}\big]\hat{w}^2(L,t)\\
&~- \epsilon_2 \left[N_2k_8-\frac{k_8}{2}-3L\big(\frac{\gamma^2}{\alpha_1}+\frac{1}{2 \beta}\big)k^2_8-\frac{1}{2}\right] \hat{p}^2(L,t)\\
&\quad-\big(L\epsilon_2 \rho -\frac{\delta_1}{2}\big) \dot{\hat{w}}^2(0,t)-\big(L\epsilon_2 \mu -\frac{\delta_2}{2}\big) \dot{\hat{p}}^2(0,t)\\
\end{array}
\end{eqnarray*}
Finally, with conditions \eqref{caps1}, \eqref{Ldot} follows.
\end{proof}

\begin{theorem}\label{thm1}
Suppose that   the constants $L, C_e,$ $\epsilon_1,$$\epsilon_2, \delta_1,\delta_2,N_1,N_2 $  fulfill the same conditions in Lemma \ref{lemm4} .
For any $\epsilon>0$ and $(v,p,\dot v,\dot p),(v_0,p_0,v_1,p_1)\in \mathcal{H},$ the energy $E(t)$  decays exponentially with the decay rate $\omega(\epsilon_1,\epsilon_2)=\min\left(\frac{\epsilon_1(1-\frac{L}{2})}{1+C_1  \epsilon_1}, \frac{\frac{2\epsilon_2}{3}}{1+C_2 \epsilon_2}\right),$ i.e.
	\begin{eqnarray}
	\label{finalresult}		&E(t)\leq \frac{p_2}{p_1} E(0)e^{-\omega (\epsilon_1,\epsilon_2) t},\qquad \forall t>0.
\end{eqnarray}
      \end{theorem}
\begin{proof}
Taking the time derivative of  $L(t)$ along with Lemma \ref{lemm4} leads to
$\dot L(t)  \le -C_e \epsilon_1 \left(1-\frac{L}{2}\right) E_e(t) - \frac{2\epsilon_2}{3}\hat{E}(t).$
Since  $ L_e (t) \le (1+ \epsilon_1 C_1) E_e (t)$ and $\hat{L}(t) \le (1+ \epsilon_2 C_2)\hat{E}(t)$ by \eqref{lemres-c1} and \eqref{lemres-c}, 
the following inequality is immediate,
\begin{eqnarray}
\begin{array}{ll}
 \dot{L}(t)
 \leq - \min \left(\frac{\epsilon_1\left(1-\frac{L}{2}\right)}{1+\epsilon_1C_1} , \frac{\frac{2\epsilon_2}{3}}{1+\epsilon_2C_2}\right) L(t).
 \end{array}
\end{eqnarray}
Now, apply the Gr\"{o}nwall's inequality  to obtain
 \begin{eqnarray}\label{punov}
     L(t)\leq \ e^{- \min \left(\frac{\epsilon_1\left(1-\frac{L}{2}\right)}{1+\epsilon_1C_1} , \frac{\frac{2\epsilon_2}{3}}{1+\epsilon_2C_2}\right)t} L(0).
 \end{eqnarray}
 Finally, by Lemma \ref{lemm2} together with \eqref{punov},
$p_1 E(t) \leq L(t) \leq e^{-\omega t} L(0) \leq p_2E(0) e^{-\omega t},$
and therefore, \eqref{finalresult} follows.
\end{proof}

{\color{black}
\section{Well-posedness}	
Letting $\vec f=[f_1,\ldots,f_8]^{\rm T}$ and $\vec g=[g_1,\ldots,g_8]^{\rm T}\in \mathrm{H}$ with  the Hilbert space  $\mathrm H:=(H^1(0,L))^4\times (L^2(0,L))^4$  equipped with the inner product
	\begin{eqnarray}
		\begin{array}{ll}
             \left<\vec f,\vec g\right>_{\mathrm H}          =\frac{1}{2}\int^L_0 \left\{ \begin{bmatrix}
				M & 0\\
				0&M
			\end{bmatrix} \begin{bmatrix}
				f_5\\
				f_6      \\
				f_7\\
				f_8
			\end{bmatrix} \cdot \begin{bmatrix}
				\bar g_5  \\
				\bar g_6       \\
				\bar g_7\\
				\bar g_8
			\end{bmatrix} \right.\\
			  \left.+ \begin{bmatrix}
				A & 0\\
				0& A
			\end{bmatrix} \begin{bmatrix}
				f_{1,x}\\
				f_{2,x}      \\
				f_{3,x}\\
				f_{4,x}
			\end{bmatrix} \cdot \begin{bmatrix}
				\bar g_{1,x}  \\
				\bar g_{2,x}       \\
				\bar g_{3,x}\\
				\bar g_{4,x}
			\end{bmatrix}  \right\} dx+ \left[\frac{k_1}{2}  f_3 \bar g_3\right.\\
			\\+ \left.\frac{k_2}{2}  f_4 \bar g_4+ \frac{k_6}{2}  f_1 \bar g_1+ \frac{k_8}{2}  f_2 \bar g_2\right] (0,t),
		\end{array}
	\end{eqnarray}
	define the system operator $\mathrm A: {\rm Dom} (\mathrm{A})\subset \mathrm H\mapsto \mathrm H$ by
	\begin{eqnarray*}
\begin{array}{ll}
\mathrm A= \begin{bmatrix}
				0 & 0 & I & 0 \\
				0 & 0 & 0& I \\
				D_x\bigotimes M^{-1} A& 0 & 0& 0  \\
				0& D_x \bigotimes M^{-1} A& 0 & 0
			\end{bmatrix},
\end{array}
\end{eqnarray*}
 		\begin{eqnarray*}
\begin{array}{ll}
{\rm Dom}(\mathrm A)=\left\{ \vec f\in \mathrm{H}: \mathrm A \vec f\in\mathrm H, \right. \\
 \begin{bmatrix}
A & 0 \\
0 & A
\end{bmatrix} \begin{bmatrix}
    f_{1,x}\\
f_{2,x}\\
 f_{3,x}\\
f_{4,x}
\end{bmatrix}(0,t) = \begin{bmatrix}
K_{13} & K_{24} \\
K_{13} & K_{24}
\end{bmatrix}\begin{bmatrix}
f_{7}\\
f_8\\
f_3\\
f_4
\end{bmatrix}(0,t), \\
\left. \begin{bmatrix}
A & 0 \\
0 & A
\end{bmatrix} \begin{bmatrix}
    f_{1,x}\\
f_{2,x}\\
 f_{3,x}\\
f_{4,x}
\end{bmatrix}(L,t) = -\begin{bmatrix}
K_{57} & K_{68} \\
0& 0
\end{bmatrix}\begin{bmatrix}
f_{5}\\
f_6\\
f_1\\
f_2
\end{bmatrix}(L,t)\right\}.
\end{array}
\end{eqnarray*}
	}
	where $D_x=\frac{\partial }{\partial x}.$
	Notice that since $M$ and $A$ are positive-definite matrices, the inner product is coercive on $\mathrm H.$ The following theorem is immediate by the L\"umer-Phillips theorem and Theorem \ref{thm1}.
\begin{lemma}
The operator $\mathrm A$ generates a $C_0-$semigroup of contractions $\{e^{{\mathrm A}t}\}$ on $\mathrm H.$
\end{lemma}

{\color{black}

\section{Proposed Model Reduction}
\label{Simu}
In alignment with the methodology outlined in \cite{Ren}, consider a given natural number $N\in\mathbb{N}$. Introduce the mesh size $h:=\frac{1}{N+1}$, and discretize the interval $[0,L]$ as follows:
$0=x_0<x_1<...<x_{N-1}<x_N<x_{N+1}=L.$
In this method, we also account for the middle points of each subinterval, denoted by $\left\{x_{i+\frac{1}{2}}\right\}_{i=0}^{N}$, where $x_{i+\frac{1}{2}}=\frac{x_{i+1}+x_i}{2}$. Define the finite average and difference operators as $u_{i+\frac{1}{2}}:=\frac{u_{i+1}+u_i}{2}, \quad \delta_x u_{i+\frac{1}{2}}:=\frac{u_{i+1}-u_i}{h}.$
It is worth noting that by considering odd-number derivatives at the in-between nodes within the uniform discretization of $[0,L]$, higher-order approximations are achieved. Furthermore, these approximations can be expressed using both the preceding and following nodes \cite{Ren}. Denoting
$(e^1, e^2, \hat w^1, \hat w^2)_j:=(e_1, e_2, \hat w,\hat p)(x_j), $ a reduced model for \eqref{systemhat}-\eqref{system-e} can be obtained as the following
\begin{eqnarray*}
	\label{ORFD2-1}
\begin{array}{ll}
			\left.\frac{M}{2}  \begin{bmatrix}
			{\ddot {\hat w}^1_{j+\frac{1}{2}}}+{\ddot {\hat w}^{1}_{j-\frac{1}{2}}}  \\
			{\ddot {\hat w}^2_{j+\frac{1}{2}}}+{\ddot {\hat w}^{2}_{j-\frac{1}{2}}}   \\
			\end{bmatrix} - \frac{A}{h} \begin{bmatrix}
			 \delta_x\left( {\hat w}^1_{j+\frac{1}{2}}- {\hat w}^1_{j-\frac{1}{2}}\right) \\
			 \delta_x \left( {\hat w}^2_{j+\frac{1}{2}}- {\hat w}^2_{j-\frac{1}{2}} \right) \\
			\end{bmatrix}\right\}_{j=1}^N=0,&\\
					\frac{M}{2}  \begin{bmatrix}
			{\ddot {\hat w}^{1}_{\frac{1}{2}}} \\
			{\ddot {\hat w}^{2}_{\frac{1}{2}} }  \\
			\end{bmatrix} - \frac{A}{h} \begin{bmatrix}
			 \delta_x {\hat w}^1_{\frac{1}{2}}\\
			 \delta_x {\hat w}^2_{\frac{1}{2}}  \\
			\end{bmatrix}  +\frac{K_{13} }{h}\begin{bmatrix}
          \dot e^1_0    \\
       \dot e^2_0
    \end{bmatrix}  +\frac{K_{24}}{h} \begin{bmatrix}
           e^1_0   \\
        e^2_0
    \end{bmatrix}=0,&\\
    \frac{M}{2}  \begin{bmatrix}
			{\ddot {\hat w}^1_{N+\frac{1}{2}}} \\
			{\ddot {\hat w}^2_{N+\frac{1}{2}}}   \\
			\end{bmatrix} - \frac{A}{h} \begin{bmatrix}
			 \delta_x {\hat w}^1_{N+\frac{1}{2}}\\
			 \delta_x {\hat w}^2_{N+\frac{1}{2}}  \\
			\end{bmatrix}  +\frac{K_{57}}{h} \begin{bmatrix}
          \dot {\hat w}^1_0    \\
       \dot {\hat w}^2_0
    \end{bmatrix}\\
    \qquad   +\frac{K_{68}}{h} \begin{bmatrix}
           \hat w^1_0   \\
        \hat w^2_0
    \end{bmatrix}=0, ~  t\in \mathbb{R}^+,&\\
 (\hat w^1, \hat w^2, \dot {\hat w}^1, \dot {\hat w}^2)_j(0)\\
 \qquad=  (w_0, p_0, w_1, p_1)(x_j), ~j=0,\ldots, N+1,\\
 	\label{ORFD2-2}\begin{array}{ll}
\left.\frac{M}{2}  \begin{bmatrix}
			{\ddot e}^{1}_{j+\frac{1}{2}}+{\ddot e}^{1}_{j-\frac{1}{2}}  \\
			{\ddot e}^{2}_{j+\frac{1}{2}}+{\ddot e}^{2}_{j-\frac{1}{2}}   \\
			\end{bmatrix} - \frac{A}{h} \begin{bmatrix}
			 \delta_x \left(e^1_{j+\frac{1}{2}}- e^1_{j-\frac{1}{2}} \right)\\
			 \delta_x \left(e^2_{j+\frac{1}{2}}- e^2_{j-\frac{1}{2}} \right) \\
			\end{bmatrix} \right\}_{j=1}^N=0,&\\
					\frac{M}{2}  \begin{bmatrix}
			{\ddot e}^{1}_{\frac{1}{2}} \\
			{\ddot e}^{2}_{\frac{1}{2}}   \\
			\end{bmatrix} - \frac{A}{h} \begin{bmatrix}
			 \delta_x e^1_{\frac{1}{2}}\\
			 \delta_x e^2_{\frac{1}{2}}  \\
			\end{bmatrix}  +\frac{K_{13}}{h} \begin{bmatrix}
          \dot e^1_0    \\
       \dot e^2_0
    \end{bmatrix}  +\frac{K_{24}}{h} \begin{bmatrix}
           e^1_0   \\
        e^2_0
    \end{bmatrix} =0,&\\
			\frac{M}{2}  \begin{bmatrix}
			{\ddot e}^{1}_{N+\frac{1}{2}} \\
			{\ddot e}^{2}_{N+\frac{1}{2}}   \\
			\end{bmatrix} - \frac{A}{h} \begin{bmatrix}
			 \delta_x e^1_{N+\frac{1}{2}}\\
			 \delta_x e^2_{N+\frac{1}{2}}  \\
			\end{bmatrix}  =0, ~  t\in \mathbb{R}^+,&\\
(e^1, e^2, \dot e^1, \dot e^2)_j(0)=0,   ~j=0,\ldots, N+1,&
\end{array}\\
\end{array}
\end{eqnarray*}
with the discretized energy $E_h(t)=\hat E_h(t)+E_{e,h}(t)$
\begin{eqnarray*}
\begin{array}{ll}
	  E_h(t):=
\frac{h}{2} \sum\limits_{j=0}^N\left\{ (I \bigotimes A)   \begin{bmatrix}
				  {\hat w}^1_{j+\frac{1}{2}} \\
				  {\hat w}^2_{j+\frac{1}{2}}      \\
				  e^1_{j+\frac{1}{2}} \\
				e^2_{j+\frac{1}{2}}       \\
			\end{bmatrix} \cdot \begin{bmatrix}
				   {\hat w}^1_{j+\frac{1}{2}}  \\
				  {\hat w}^2_{j+\frac{1}{2}}        \\
				    e^1_{j+\frac{1}{2}} \\
				e^2_{j+\frac{1}{2}}       \\
			\end{bmatrix}  \right.\\
\left.			
	(I \bigotimes \frac{M}{2}) \begin{bmatrix}
				{\dot {\hat w}}^1_{j+\frac{1}{2}} +{\dot {\hat w}}^1_{j-\frac{1}{2}}  \\
				{\dot {\hat w}}^2_{j+\frac{1}{2}}  +{\dot {\hat w}}^2_{j-\frac{1}{2}}       \\
				{\dot {e}}^1_{j+\frac{1}{2}}+{\dot {e}}^1_{j-\frac{1}{2}}  \\
				{\dot {e}}^2_{j+\frac{1}{2}} +   	{\dot {e}}^2_{j-\frac{1}{2}}    \\
			\end{bmatrix} \cdot \begin{bmatrix}
					{\dot {\hat w}}^1_{j+\frac{1}{2}} +{\dot {\hat w}}^1_{j-\frac{1}{2}}  \\
				{\dot {\hat w}}^2_{j+\frac{1}{2}}  +{\dot {\hat w}}^2_{j-\frac{1}{2}}       \\
				{\dot {e}}^1_{j+\frac{1}{2}}+{\dot {e}}^1_{j-\frac{1}{2}}  \\
				{\dot {e}}^2_{j+\frac{1}{2}} +   	{\dot {e}}^2_{j-\frac{1}{2}}    \\
			\end{bmatrix}		 \right\} dx\\
			 +\frac{1}{2} \left( k_6 |{\dot {\hat w}}^1_{N+\frac{1}{2}}|^2+k_8 |{\dot {\hat w}}^2_{N+\frac{1}{2}}|^2 +k_2 |{\dot {e}}^1_{\frac{1}{2}}|^2 + k_4|{\dot {e}}^2_{\frac{1}{2}}|^2\right),
\end{array}\end{eqnarray*}
where $\bigotimes$ is the Kronecker matrix product.  In the algorithm implementation, assuming no observation error,  we set $T=5$ seconds, $L=1$ meter, and $N=30$ (or $h=1/31)$. Material parameters are $\rho=1$ kg/m³, $\mu=10^{-1}$ H/m, $\alpha=1$ N/m², $\gamma=10^{-2}$ C/m³, and $\beta=3$ m/F.

Initial conditions $(w_0,w_1,p_0,p_1)(x)$ are chosen intentionally a  high-frequency-type: $\frac{10^{-3}}{25} \sum\limits_{k=5}^{30} \cos\left(\frac{k\pi x}{L}\right)$. Controllers are set to $k_1\to 10^{-7}$, $k_2\to 10^{-8}$, $k_3\to 10^{-7}$, $k_4\to 3 \times 10^{-6}$, $k_5\to 60$, $k_6\to 2\times 10^{-2}$, $k_7\to 10$, $k_8\to 4\times 10^{-2}$. In Fig. \ref{3dplots}, solutions decay to zero exponentially fast.
\begin{figure}[htb!]
\includegraphics[width=1.7in]{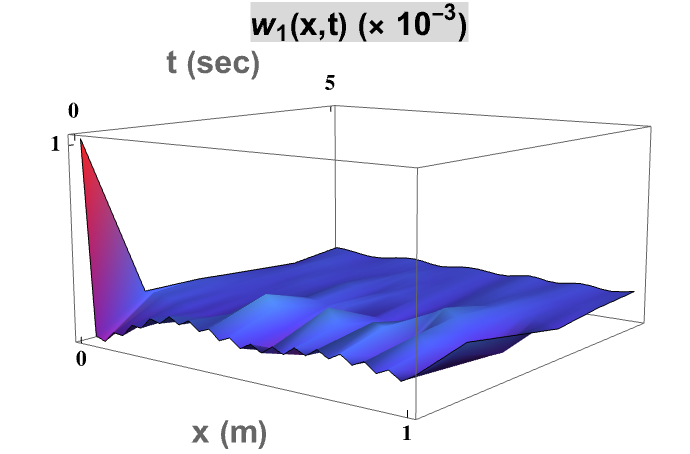} \includegraphics[width=1.7in]{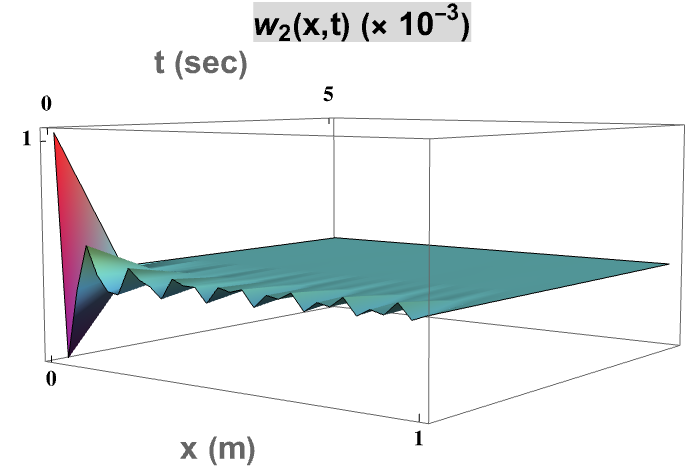}
\includegraphics[width=1.7in]{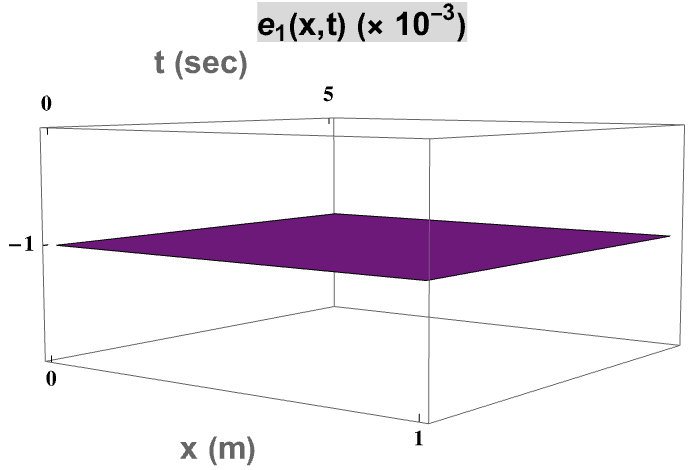} \includegraphics[width=1.7in]{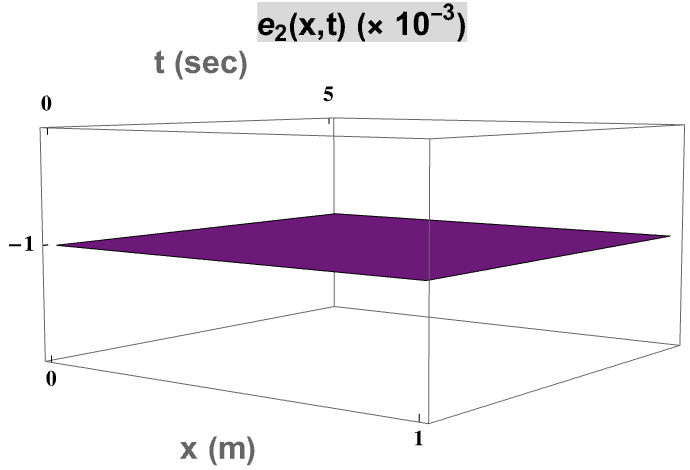}
\caption{\footnotesize Solutions are observed to decay to zero exponentially fast.}
	\label{3dplots}
\end{figure}

\section{Conclusion\& Future Work}
The comprehensive analysis presented here is relevant for the control of multi-layer stack designs \cite{IFAC}. Moving forward, our immediate focus is on analyzing the solutions derived from the order-reduction algorithm, as outlined in Section \ref{Simu}, following the methodology presented in \cite{Ren}. The proposed algorithm demonstrates the feasibility of employing a discrete Lyapunov approach to establish the exponential stability of the order-reduced system,  uniformly as $h\to 0$ \cite{U2}.

\end{document}